\documentclass[12pt]{article}
\usepackage[centertags]{amsmath}
\usepackage{amsfonts}
\usepackage{amssymb}
\usepackage{amsthm}
\usepackage{newlfont}
\usepackage[latin1]{inputenc}
\def\RR{{\mathrm I}\!{\mathrm R}}

\def\pe{\textbf{P} }
\def\besov{||.||_{BE^{\beta}}}
\newtheorem{fed}{Definition}[section]
\newtheorem{teo}[fed]{Theorem}
\newtheorem{lem}[fed]{Lemma}
\newtheorem{coro}[fed]{Corollary}

\newtheorem{defn}[fed]{Definition}

\begin{document}
\title{Stability of bounded global solutions for Navier-Stokes equations
 }



\author{ Oscar A. Barraza \quad Claudia B. Ruscitti}
\date{Mayo, 2007.}
\maketitle \vskip-2cm
\ \\
\begin{center}
{\it Departamento de Matem\'atica, Facultad de Ciencias Exactas,
 Universidad Nacional de La Plata, C.C. 172, (1900) La Plata,
 Argentina.\hfill
 e-mail: oscar@mate.unlp.edu.ar, claudia@mate.unlp.edu.ar}
\end{center}
\vskip3cm




\begin{abstract}

In this paper some kind of asymptotic behavior of the solutions
for the Navier-Stokes system on $\RR ^n $ in abstract Banach
spaces is studied under the existence of global in time solutions.
The asymptotic stability of the zero solution is also shown.

\end{abstract}

\vfill \eject

\section{Introduction}

\

The Navier -Stokes equations describing the  motion of an
incompressible fluid in $\RR^n,\ n \ge 2$, without external forces
are written as follows

\begin{eqnarray}\label{fuerte}
u_t -\triangle u +(u. \nabla) u + \nabla p  = 0 &\nonumber \\
 \nabla . u =0, &\\
u(0)=u_0.& \nonumber
\end{eqnarray}

Here  $u(x,t)=(u_1(x,t),...,u_n(x,t))$ is the unknown velocity of
the fluid and $p = p(x,t)$ is the unknown pressure at the point $x
\in \RR^n$ and time $t\ge 0$.

Let $u_0(x)$ be an initial condition that verifies $\nabla . u_0
=0$.

As usual, $\textbf{P}$ denotes the projection from
$(L^2(\RR^n))^n$ onto the subspace $\textbf{P}(L^2(\RR^n))^n=\{ f
\in L^2: \nabla . f = 0\}$ of solenoidal vector fields

\[\textbf{P} (u_1,...,u_n) = (u_1- R_1 \sigma,..., u_n-R_n
\sigma),\]

where $R_j$ are the Riesz transforms, with symbols $\xi_j/
|\xi_j|$ and $\sigma= R_1u_1 +...+ R_n u_n$.

It is well known that $\textbf{P}$ can be extended to a continuous
operator on $L^p(\RR^n) $ to $\pe(L^p(\RR^n)) =\{ f \in L^p
(\RR^n) : \nabla . f =0 \}$. ( c fr \cite{f-m})

Using $\textbf{P}$ the equations (\ref{fuerte}) can be rewritten
as

\begin{eqnarray}\label{proyectada}
u_t &=& \triangle u -\textbf{P} \nabla (u\otimes u)\nonumber\\
 \nabla\  . \ u &=&0 \\
u(0)&=&u_0. \nonumber
\end{eqnarray}

 $(u.\nabla)u$ can be replaced by $\nabla (u\otimes u)$, because $\nabla . u=0$. \

The heat semigroup $S(t)$ is given by the convolution with the
Weierstrass kernel, or heat kernel, $\displaystyle{S(t)= e^{t
\triangle}= \frac{1}{(4 \pi t)^{n/2}}\ e^{-|x|^2/{4t}}\ *}$; then
the problem (\ref{proyectada}) can be written under the following
integral form

\begin{equation}
u(t) = S(t)u_0 - \int_0^t \textbf{P} \nabla S(t-  \tau) (u\otimes
u)(\tau)d \tau \label{mild}
\end{equation}

where the integrals are understood in the Bochner's sense. The
solutions of this integral equation are called mild solutions.

 Then, a solution of (\ref{fuerte}) or (\ref{proyectada})
will be interpreted as an $E$ valued mapping defined in $[0,
\infty)$, for an appropriate Banach space $E$.

In this article, by solutions of the Navier-Stokes system we mean
solutions of type (\ref{mild}).

 The  goal of this paper is to show a stability
result for global in time solutions for (\ref{fuerte}). In theorem
2 the difference of two given
 bounded solutions in an abstract Banach space
goes to zero with the precise rate of decay. As an immediate
consequence the asymptotic stability of the zero solution is
obtained.

 In section two we present the mathematical setting to
build the solutions mentioned above. Section three is dedicated to
state and prove our main theorem.

\section{Abstract Banach spaces}

In this section we retrieve all the needed definitions to
construct a Banach space adequate to the Navier-Stokes system.
These ideas were introduced in  \cite{c-k}, \cite{k} and even
\cite{lr}.

\begin{defn}
 A Banach space $(E, ||.||_E)$ is said to be functional and translation
invariant if the following three conditions are satisfied:

(i) $S \subset E \subset S' $ and  both inclusions are continuous,

(ii) for every $f \in E, \ \tau : \RR^n \to E$ defined by $\tau_y
f(x) = f(x+y)$ is measurable in the sense of Bochner with respect
to the Lebesgue measure on $\RR^n$,


(iii) the norm $||.||_E$ on $E$ is invariant  translation

$$\forall f \in E, \ y \in \RR^n\ ||\tau_y f||_E
= ||f||_E .$$

\end{defn}

\begin{defn}

 We call the space $(E, ||.||_E)$ adequate to the problem
(\ref{fuerte}) if

(i) $(E, ||.||_E)$ is a functional and translation invariant
Banach space.

(ii) $\forall f, g \in E$ the product $f \otimes g$ is well
defined as a tempered distribution. Moreover, there exists $T_0
> 0$ and a positive function $w\in L^1 (0, T_0)$ such that

$$ ||\pe \nabla S(\tau)(f\otimes g)||_E \le w(\tau) ||f||_E
||g||_E, \ \forall f, g\in E, \tau \in (0, T_0).$$

\end{defn}

Some examples of Banach spaces adequate to (\ref{fuerte}) are the
subspaces of free divergence functions of the spaces $L^p,\
L^p_w,\ L^{p,q}, \ M^p_q$ with $p>n$ and $q\ge 1$.

In this paper, the norms of the Banach spaces have additional
scaling properties. Let
 $f: \RR^n \to \RR^n$ be, the rescaled function  $f_{\lambda}$ is defined by

$$f_{\lambda}(x)= f(\lambda x), \ \lambda >0.$$

This definition is extended for all $f \in S'$ in the usual way.

\begin{defn}

 Let $(E, ||.||_E)$ be a Banach space which can be
imbedded continuously in $S'$. The norm $||.||_E$ is said to have
a scaling degree equal to k if

\

$||f_{\lambda}||_ E = \lambda ^k ||f||_E, \ \forall f \in E $ such
that $ f_{\lambda} \in E$ and $\forall \lambda > 0.$

\end{defn}

{\bf Remark:} The usual norms of the spaces $L^p,\ L^p_w,\
L^{p,q}, \ M^p_q$ have scaling degree equal to $-n/p$. Moreover,
the standard norm in  the homogeneous Sobolev space $\dot{H}^s =\{
f \in S': |\xi|^s \hat{f}(\xi) \in L^2\}$ has scaling degree equal
to $s-n/2$.

\

A Banach space $E$ endowed with a norm with scaling degree equal
to $-n/p,\ p>n$ will be denoted by $E_p$.



 Let $E\subset S'$ be a Banach space, in \cite{k} was introduced the space of
 distributions $BE^{\beta}$
 which is an `` homogeneous  Besov space modelled on $E$''.

 \begin{defn}
 Let $\beta \ge 0$. Given a Banach space $E$ continuously imbedded in $S'$, define

 $$BE^{\beta}= \{f \in S': ||f||_{BE^\beta} = \sup_{t>0}
 t^{\beta/2}||S(t)f||_E < \infty\}.$$
\end{defn}

 For instance, for the Banach space $E= L^p$, the norm $  ||.||_{BE^{\beta}}$ is
 equivalent to the standard norm of the homogeneous Besov spaces ${\dot B}^{-\beta}_{p,
 \infty}$.

It is straightforward to see that if $E$ has a norm with scaling
degree equal to $k$, then $\besov$ has scaling degree equal to
$k-\beta$. Indeed, given any $f\in S'$ and $\lambda >0$,
$$S(t)f_{\lambda}= {(S({\lambda t}^2)f)}_{\lambda}$$
then,
$$||f_\lambda||_{BE^{\beta}}= \sup_{t>0}t^{\beta /
2}||S(t)f_{\lambda}||_E= {\lambda}^{k-\beta} \sup_{\lambda^2 t >0}
{(\lambda^ 2 t)}^{\beta/2} ||S(\lambda ^2 t) f||_E = \lambda ^{k -
\beta} ||f||_ {BE^\beta}.$$

 It is well known that for the case of a Banach space $E_p$  such that for some
 $q \in [1, \infty],\  e^{t\triangle} : E_p \to
 L^q$ is a bounded operator for every $t>0$, then
 $(BE_p^{\beta}, \ ||.||_{BE_p^\beta})$ is a Banach space.

These definitions allow us to construct global in time solutions
(for small initials conditions) in the space $\chi \equiv {\cal
C}([0, \infty), BE_p^{\beta})$ consisting of measurable and
essentially bounded functions $u: [0, \infty) \to BE_p^{\beta}$
such that $u(t) \to u(0)$ as $t\searrow 0$ in the topology de
$S'$.


\section{Global in time solutions to the Navier-Stokes equations}

In this section the main result is exposed. For this purpose, we
recall the following  two statements. In the first lemma the
boundedness property of the operator $\textbf{P} \nabla S(t)
(u\otimes v)= e^{t\triangle}B(u,v)$, where $B(u,v)=\pe \nabla (u
\otimes v)$, is shown.

Theorem 1 gives the conditions for the existence and unicity of
global solutions. Although this result is independent of our main
statement (Theorem 2), we decided to include it here just to keep
in mind that the set of global solutions of system (\ref{fuerte})
is neither empty nor single (with zero as the unique global
solution).

\begin{lem}\cite{k}
Assume that the Banach space $E_p$ is adequate to the problem
(\ref{fuerte}) and has a norm with scaling degree equal to $-n/p$.
Then

\begin{enumerate}
\item There exists a constant $C_1>0$ independent of $t,\ u,\ v$,
such that

$$||e^{t\triangle}B(u,v)||_{E_p} \le C_1 t^{-(1+n/p)/2}
||u||_{E_p} . ||v||_{E_p}$$

for all  $u,\ v \in E_p\  and\  t>0$.

\item Let $0\le \beta\le 1+ n/p$.  There exists a constant $C_2
>0$, independent de $t,\ u,\ v$ such that

$$||e^{t\triangle}B(u,v)||_{B E_p^{\beta}} \le  C_2 t^{(\beta-1-n/p)/2}
||u||_{E_p} . ||v||_{E_p}$$

\end{enumerate}
\end{lem}

In the following theorem  global in time solutions are built for
the problem (\ref{fuerte}) in the space $BE_p^{\beta}$ considering
small initial conditions. The proof of this result is obtained by
handing standard tools introduced in \cite{ka} and used by other
authors. See, for instance, \cite{b}, \cite{c}, \cite{c-k}, or
\cite{lr}.

\begin{teo}\cite{k} Fix $p>n,\ n>0$. Denote $\beta = 1-n/p$. Let $E_p$ be a
Banach space satisfying

\
 (i) $E_p$ is an adequate space  to the problem (\ref{fuerte});

 (ii) the norm $||.||_{E_p}$ has scaling degree equal to $-n/p$;

 (iii) there exists  $q>0$ such that the operator $e^{t\triangle}: E_p \to
 L^q$ is bounded for every $t>0$.

 Then, there exists  $\epsilon >0 $ such that for every $u_0 \in BE_p^{\beta}$
 satisfying $||u_0||_{BE_p^{\beta}} < \epsilon$ there exists a solution of
 (\ref{fuerte}) for all $t>0$ in the space

 $$\chi\equiv {\cal C}([0, \infty), BE_p^{\beta}) \cap \{ u: (0, \infty) \to E_p :
 \sup_{t>0} t^{\beta/2}||u(t)||_{E_p} < \infty\}.$$

This is the unique solution satisfying the condition $\sup_{t>0}
t^{\beta/2}||u(t)||_{E_p}\le 2 \epsilon$.

\end{teo}
Now, we are able to state our main result.

\begin{teo} Fix $p>n,\ n>0$. Denote $\beta = 1-n/p$. Let $E_p$ be a
Banach space satisfying the same hypothesis (i), (ii),(iii) as in
the previous theorem. Let  $u,\ v$ be two global in time solutions
of (\ref{fuerte}) in the space $\chi$ corresponding to the initial
conditions $u_0,\ v_0 \in BE_p^\beta$ respectively. Suppose that
$||u||_{\chi} \le M$ and $||v||_{\chi} \le M$, for some constant
$M>0$.

Let $w_0= u_0 -v_0$ and $w(t) = u(t) -v(t)$.

Then
$$||w (t)||_{E_p} \le (||w_0||_{BE_p^{\beta}} + 4CM^2) t^{-\beta/2}, \
t>0.$$

That is, $||u(t) -v(t)||_{E_p}$ goes asymptotically  to 0 like
$t^{-\beta/2}$.

\

\end{teo}

{\bf Remark:} This theorem plays a fundamental role in providing a
way to compare two bounded global in time solutions  for the
system (\ref{fuerte}), corresponding to any initial conditions,
even if the norms of these initial conditions are not small
enough.

{\it Proof}

Let the space

 $$\chi\equiv {\cal C}([0, \infty), BE_p^{\beta}) \cap \{ u: (0, \infty) \to E_p :
 \sup_{t>0} t^{\beta/2}||u(t)||_{E_p} < \infty\}$$

be endowed  with the norm

$$||u||_{\chi} = \max \{ \sup_{t>0} ||u(t)||_{BE_p^{\beta}} ,\
\sup_{t>0} t^{\beta/2}||u(t)||_{E_p} \}$$

 Let $u,\ v$ two global in time solutions of (\ref{fuerte})
 corresponding to the initial data $u_0,\ v_0 \in BE_p^{\beta}$,
 respectively.

 We know that
\begin{eqnarray*}
u(t) &= S(t)u_0 - \displaystyle{\int_0^t} \textbf{P} \nabla S(t-
\tau) (u\otimes
u)(\tau)d \tau \\
&= e^{t\triangle} u_0 - \displaystyle{\int_0^t} e^{(t-\tau)
\triangle} B(u,u)
d\tau.\\
\end{eqnarray*}

Analogously,
$$v(t)= e^{t\triangle} v_0 - \displaystyle{\int_0^t} e^{(t-\tau) \triangle} B(v,v)
d\tau$$ and hence,
$$u(t)- v(t)= e^{t\triangle} (u_0 -v_0) \ -\ \displaystyle{\int_0^t} e^{(t-\tau)
\triangle}[B(u,u) -B(v,v)] d \tau.$$

From this, we have

\begin{tabular}{ll}
$t^{\beta/2} ||u(t) - v(t)||_ {E_p}$& $= t^{\beta/2} ||w(t)||_
{E_p} $ \\
\\
 & $\le t^{\beta/2}
||e^{t\triangle}(u_0 - v_0)||_{E_p} $ \\
\\
&\ \ $ + t^{\beta/2} \displaystyle{\int_0^t} ||e^{(t-\tau)
\triangle} [B(u,u - v) +
B(u-v, v)]||_{E_p} d\tau $\\
 \\

& $ \le||w_0||_{BE_p^{\beta}}$ \\
\\
& \ \ $ + t^{\beta/2} \displaystyle{\int_0^t} [||e^{(t-\tau)
\triangle} B(u,u - v)||_{E_p} + ||e^ {(t-\tau)
\triangle}B(u-v, v)||_{E_p}] d\tau .$\\
 \\
\end{tabular}

Applying lemma 1 we can conclude

\begin{tabular}{ll}
$t^{\beta/2} ||u(t) - v(t)||_ {E_p}$&
$\le ||w_0||_{BE_p^{\beta}} $\\
\\
& \ \ $ + t^{\beta/2} \displaystyle{\int_0^t} C (t-\tau)^{-(1+
n/p)/2}
[||u(\tau)||_{E_p} ||(u - v)(\tau)||_{E_p}  $ \\
\\
&\ \ $ +||v(\tau)||_{E_p} ||(u - v)(\tau)||_{E_p}] d\tau $\\
 \\
&$ \le ||w_0||_{BE_p^{\beta}}  $\\
\\
&\ \ $+ C (||u||_{\chi} + ||v||_{\chi}) t^{\beta/2}
\displaystyle{\int_0^t} (t-\tau)^{-(1+ n/p)/2} \tau^{-\beta}(
\tau ^{\beta /2} ||(u - v)(\tau)||_{E_p} ) d \tau $\\
 \\
&$\le ||w_0||_{BE_p^{\beta}}  + C (||u||_{\chi} + ||v||_{\chi})
(\sup _{\tau >0} \tau ^{\beta /2} ||(u - v)(\tau)||_{E_p}) $\\
 \\
&$\le ||w_0||_{BE_p^{\beta}} +  C (2M)^2.$\\
\end{tabular}

Therefore,
$$  ||w(t)||_
{E_p} =||u(t) - v(t)||_ {E_p} \le (||w_0||_{BE_p^{\beta}} + 4 C M
^2) t^{ -\beta/2}.$$

\

\begin{coro} Under the hypotheses of Theorem 2 the zero solution
is asymptotically stable. More precisely, let $u$ be a bounded
global solution of (\ref{fuerte}) with initial condition $u_0$,
the norm $||u||_{E_p}$ tends to zero as $t^{-\beta/2}$.
\end{coro}

{\it Proof}

It is immediately to show this result after taking $v = 0$ in the
previous theorem.

{\bf Remark:} This property of the zero solution holds in almost
all the known solutions spaces for the Navier-Stokes equations.

\begin{coro} Let the assumptions from the  Theorem 2
hold. Let $u,\ v$ be two solutions of (\ref{fuerte}) corresponding
to the initial data $u_0,\ v_0 \in BE_p^\beta$ respectively.
Suppose that

$$\lim_{t \to \infty} t^{\beta/2} ||e^{t\triangle} (u_0 -
v_0)||_{E_p} = 0.$$ Then,
$$\lim_{t \to \infty} t^{\beta/2} ||u(.,t) -v(.,t)||_{E_p} = 0$$

provided that $M > 0$ is sufficiently small.
\end{coro}

Before proving this statement, let us recall the following
technical lemma.

\begin{lem}\cite{b} Let $w \in L^1 (0,1)$ and $\int_0^1 w(x)dx < 1$. Assume
that $f$ and $g$ are two nonnegative, bounded measurable functions
such that

$$ f(t)\le g(t)+ \int_0^1w(\tau) f(\tau t) d \tau.$$

Then $\lim _{t \to \infty} g(t)= 0$ implies $\lim _{t \to \infty}
f(t)= 0$.
\end{lem}

{\it Proof of Corollary 2}

As in Theorem 2,

\begin{tabular}{ll}
$t^{\beta/2} ||u(t) - v(t)||_ {E_p}$& $\le t^{\beta/2}
||e^{t\triangle}(u_0 - v_0)||_{E_p} $\\
\\
& \ \ $ + t^{\beta/2} \displaystyle{\int_0^t} C (t-\tau)^{-(1+
n/p)/2}
[||u(\tau)||_{E_p} ||(u - v)(\tau)||_{E_p}  $ \\
\\
&\ \ $ +||v(\tau)||_{E_p} ||(u - v)(\tau)||_{E_p}] d\tau .$\\
 \\
\end{tabular}

Since the solutions $u,\ v$ are bounded with bound $M$, we have

\begin{tabular}{ll}
$t^{\beta/2} ||u(t) - v(t)||_ {E_p}$
&$\le t^{\beta/2}||e^{t\triangle}(u_0 - v_0)||_{E_p} $\\
\\
& \ \ $+ 2MC  t^{\beta/2} \displaystyle{\int_0^t} (t-\tau)^{-(1+
n/p)/2} \tau^{-\beta}(
\tau ^{\beta /2} ||(u - v)(\tau)||_{E_p} ) d \tau .$\\
 \\
\end{tabular}

After a change of variables, it is possible to write

\begin{tabular}{ll}
$t^{\beta/2} ||u(t) - v(t)||_ {E_p}$
&$\le t^{\beta/2} ||e^{t\triangle}(u_0 - v_0)||_{E_p}$\\
\\
& \ \ $ + 2CM \displaystyle{\int_0^1} (1-s)^{- (1+n/p)/2}
s^{-\beta} (ts)^{\beta/2} ||(u-v)(ts)||_{E_p}ds.$

\end{tabular}

Putting $f(t) =  t^{\beta/2} ||u(t) - v(t)||_ {E_p}$ and using
that $(1-s)^{- (1+n/p)/2} s^{-\beta}  \in L^1 (0,1) $ we may apply
Lemma 3 to obtain $t^{\beta/2} ||u(t) - v(t)||_ {E_p} \to 0 $ as
$t \to \infty$ for sufficiently small $M >0$.


\end{document}